\documentclass[11pt]{article}
\usepackage{amstext,amssymb,amsmath,amsbsy}

\usepackage{tikz}

\usepackage{hyperref}
\usepackage{amscd}
\usepackage{amsfonts, dsfont}
\usepackage{indentfirst}
\usepackage{verbatim}
\usepackage{amsmath}
\usepackage{amsthm}
\usepackage{enumerate}
\usepackage{graphicx}
\usepackage{subfig}
\usepackage{color}
\usepackage[OT1]{fontenc}
\usepackage[latin1]{inputenc}
\usepackage[english]{babel}
\usepackage{amssymb}
\newtheorem{theorem}{Theorem}
\newtheorem{lemma}{Lemma}

\newtheorem{proposition}{Proposition}

\newtheorem{remark}{Remark}

\numberwithin{equation}{section}

\newcommand{\proofend}{\hfill $\Box$ }

\newcommand{\bH}{{\bf H}}
\newcommand{\cE}{{\cal E}}
\newcommand{\cH}{{\cal H}}

\newcommand{\oE}{E^{(1)}}
\newcommand{\oH}{H^{(1)}}
\newcommand{\tE}{E^{(2)}}
\newcommand{\tH}{H^{(2)}}

\newcommand{\ocE}{\widetilde E}
\newcommand{\ocH}{ \widetilde H }

\newcommand{\obE}{{\bf E}^{(1)}}
\newcommand{\obH}{{\bf H}^{(1)}}

\newcommand{\Beq}{\begin{equation}\left\{\begin{array}{rcll}}
\newcommand{\Eeq}[1]{\end{array}\right.\label{#1}\end{equation}}
\newcommand{\Eeqq}{\end{array}\right.\nonumber\end{equation}}
\newcommand{\beq}{\begin{equation}}
\newcommand{\eeq}[1]{\label{#1}\end{equation}}

\newcommand{\pB}{\partial B}

\newcommand{\supp}{\operatorname{supp}}

\newcommand{\dive}{\operatorname{div}}
\newcommand{\curl}{\operatorname{curl}}

\newcommand{\eps}{\varepsilon}

\newcommand{\loc}{_{loc}}

\newcommand{\mC}{\mathbb{C}}

\newcommand{\mN}{\mathbb{N}}

\newcommand{\mR}{\mathbb{R}}

\newcommand{\teps}{{\widetilde \eps}}
\newcommand{\tmu}{{\widetilde \mu}}
\newcommand{\sign}{\mbox{sign }}

\title{Cloaking using complementary media for electromagnetic waves}

\author{Hoai-Minh Nguyen \footnote{EPFL SB MATHAA CAMA, Station 8,  CH-1015 Lausanne, hoai-minh.nguyen@epfl.ch}  %\footnote{The research is partially supported by NSF grant DMS-1201370 and by the Alfred P. Sloan Foundation.}
}

\begin{document}
%\date{}

\maketitle 

\begin{abstract}
Negative index materials are artificial structures whose refractive index has negative value over some frequency range. The study of these materials has attracted a lot of attention in the scientific community not only because of  their many potential interesting applications but also because of challenges in understanding their intriguing properties due to the sign-changing coefficients in  equations describing their properties. In this paper, we establish cloaking using complementary media for electromagnetic waves. This confirms and extends the suggestions in two dimensions  of Lai et al. in \cite{LaiChenZhangChanComplementary} for the full Maxwell equations. The analysis is based on the reflecting and removing localized singularity techniques,  three-sphere inequalities, and the fact that  the Maxwell equations can be reduced to a weakly coupled second order elliptic equations. 

\end{abstract}



\section{Introduction}

Negative index materials (NIMs) are artificial structures whose refractive index has negative value over some frequency range. These materials were  investigated theoretically by Veselago in \cite{Veselago}.  The  existence of such materials was confirmed by Shelby, Smith, and Schultz in \cite{ShelbySmithSchultz}.  The study of NIMs has attracted a lot of attention in the scientific community not only because of  their many potential interesting applications but also because of challenges in understanding  intriguing properties of these materials.

One of  the interesting applications of NIMs is cloaking using complementary media, which 
was inspired by the concept of complementary media, see \cite{PendryRamakrishna, LaiChenZhangChanComplementary, Ng-Complementary, Ng-Superlensing-Maxwell}.  
Cloaking using complementary media was proposed and studied numerically by Lai et al. in \cite{LaiChenZhangChanComplementary} in two dimensions. The idea of this cloaking technique is to cancel the light effect of an object using its complementary media.  
Cloaking using  complementary media was mathematically  established by Nguyen in \cite{Ng-Negative-Cloaking} in the quasistatic regime.  To this end, we introduced the removing localized singularity technique and used  a standard three-sphere inequality. The method used in \cite{Ng-Negative-Cloaking} also works for the Helmholtz equation. Nevertheless, it requires small size of the cloaked region for large frequency due to the use of the (standard) three-sphere inequality. In \cite{MinhLoc2}, Nguyen and (H. L.) Nguyen gave a mathematical proof of cloaking using complementary media in the finite frequency regime for acoustic waves without imposing any condition on the size of the cloaked region. To successfully apply the approach in \cite{Ng-Negative-Cloaking}, we established a new three-sphere inequality for the Helmholtz equations which holds for arbitrary radii. 
%A modification of the cloaking setting  to obtain illusion optics with roots from \cite{Ng-Superlensing} is also discussed there.  

Another cloaking object technique  using NIMs is cloaking an object via anomalous localized resonance technique. This was suggested and studied by Nguyen in \cite{Ng-CALR-object}. Concerning this technique, an object is cloaked by the complementary property (or more precisely by the doubly complementary property) of the medium. 
This cloaking technique is inspired by the work of Milton and Nicorovici in \cite{MiltonNicorovici}. In their work, they discovered  cloaking a source  via anomalous localized resonance for constant radial plasmonic structures in the two-dimensional quasistatic regime (see \cite{BouchitteSchweizer10,AmmariCiraoloKangLeeMilton2, KohnLu, Ng-CALR-CRAS, Ng-CALR, Ng-CALR-frequency} and the references therein for recent results in this direction).  Another interesting application of NIMs is superlensing,  i.e., the possibility to beat the Rayleigh diffraction limit: no constraint between the size of the object and the wavelength is imposed, see \cite{Ng-Superlensing, Ng-Superlensing-Maxwell} and references therein. 
%It is worthy noting that a lens can act like a cloak and conversely,  see \cite{Ng-CALR-object}.

Two difficulties in the study of cloaking using complementary are as follows. Firstly, the problem is unstable. This can be explained by the fact that the equations describing the phenomena  have sign-changing coefficients; hence the ellipticity and the compactness are lost in general. Indeed, this is the case for complementary media, see \cite{Ng-WP}.   Secondly, localized resonance might appear, i.e., the field explodes in some regions and remains bounded in some others.  
It is worthy noting that the character of resonance associated with NIMs is quite complex; localized resonance and complete resonance can occur, see \cite{MinhLoc1}.

In this paper, we study cloaking using complementary media for electromagnetic waves (Theorem~\ref{thm1}).
This is a natural continuation of \cite{Ng-Negative-Cloaking, MinhLoc2} where the acoustic setting was investigated. Moving from the acoustic setting to electromagnetic one makes the analysis more delicate due to the lack of the  information of the whole gradient in the Maxwell equations. Besides extending the removing localized singualirty technique in \cite{Ng-Negative-Cloaking, Ng-Superlensing} with a simplification in \cite{Ng-CALR-frequency}, and using the reflecting technique in \cite{Ng-Superlensing-Maxwell} with roots in \cite{Ng-Complementary}, the analysis also involves recent stability result for Maxwell equations in \cite{Ng-Superlensing-Maxwell}, the weakly coupled second order elliptic equations property of Maxwell equations, and  new three-sphere inequalities. Using the concept of reflecting complementary in \cite{Ng-Superlensing-Maxwell},  one can establish cloaking using complementary media for the Maxwell equations for a general class of schemes (Proposition~\ref{pro1} in Section~\ref{sect-further-discussion}).

Let us now describe in details a scheme to  cloak an arbitrary object using complementary media for the Maxwell equations.  A more general class of schemes is considered in Section~\ref{sect-further-discussion}. 
Let  $B_r$ denote the ball centered at the origin and of radius $r$ in $\mR^3$ unless specified otherwise and let $\langle \cdot, \cdot, \rangle$ denote the Euclidean  scalar product in $\mR^3$.  Assume that the cloaked region is the annulus  $B_{2r_2} \setminus B_{r_2}$ in $\mR^3$ for some $r_{2}> 0$ in which the medium is characterized by a pair of two matrix-valued functions $(\eps_O, \mu_O)$ of  the permittivity $\eps_O$ and the permeability $\mu_O$ of the region. The assumption on the cloaked region by all means imposes no restriction  since any bounded set is a subset of such a region provided that the radius and the origin are appropriately chosen.  We assume that $\eps_O$ and $\mu_O$ are uniformly elliptic,  i.e., 
\begin{equation}\label{ellipticity}
\frac{1}{\Lambda} |\xi|^2 \le \big\langle \eps_O(x) \xi, \xi \big\rangle \le \Lambda |\xi|^2 \mbox{  and  } \frac{1}{\Lambda} |\xi|^2 \le \big\langle \mu_O(x) \xi, \xi \big\rangle \le \Lambda |\xi|^2  \quad \forall \, \xi \in \mR^3, \mbox{ a.e. } x \in B_{r_2 } \setminus B_{r_1}.  
\end{equation}
In this paper, we use a  scheme in \cite{Ng-Negative-Cloaking} with roots in the work of Lai et al. \cite{LaiChenZhangChanComplementary}. Following \cite{Ng-Negative-Cloaking}, the cloak contains two parts. The first one, in $B_{r_2} \setminus  B_{r_1}$, makes use of complementary media to cancel the effect of the cloaked region and the second one, in $B_{r_1}$, is to fill the space which ``disappears" from the cancellation by the homogeneous medium. Concerning the first part, 
instead of $B_{2r_2} \setminus B_{r_2}$, we consider $B_{r_3} \setminus B_{r_2}$ for some $r_3 > 0$ as the cloaked region in which the medium is given by 
\begin{equation*}
\big(  \teps_O, \tmu_O \big) = \left\{ \begin{array}{cl} \big(\eps_O, \mu_O\big) & \mbox{ in } B_{2 r_2} \setminus B_{r_2}, \\[6pt]
\big(I, I \big) & \mbox{ in } B_{r_3} \setminus B_{2 r_2}. 
\end{array} \right. 
\end{equation*} 
(This extension first seems to be technical but is later proved to be necessary in the acoustic setting, see \cite{Ng-CALR-object}.)  The (reflecting) complementary medium  in $B_{r_2} \setminus B_{r_1}$ is then given by 
\begin{equation}\label{first-layer}
\big(F^{-1}_*\teps_O, F^{-1}_*\tmu_O \big), 
\end{equation}
where $F: B_{r_2} \setminus \bar B_{r_1} \to B_{r_3} \setminus \bar B_{r_2}$ is the Kelvin transform with respect to $\partial B_{r_2}$, i.e., 
\begin{equation}\label{def-F}
F(x) = \frac{r_2^2}{|x|^2} x. 
\end{equation}
Here
\begin{equation}\label{def-T*}
{\cal T}_* a (y) = \frac{\nabla {\cal T}  (x)  a(x) \nabla  {\cal T} ^{T}(x)}{J(x)}, 
\end{equation}
where $x ={\cal T}^{-1}(y)$ and $J(x) = \det \nabla  {\cal T}(x)$  for a diffeomorphism ${\cal T}$. It follows that 
\begin{equation}\label{cond-r}
 r_1 = r_2^2/ r_3.   
\end{equation}
Note that in the definition of $T_*$ given in \eqref{def-T*}, $J(x): = \det \nabla T(x) $ {\bf not} $|\det \nabla T(x)|$ as often used in the acoustic setting. This convention is very suitable for the  electromagnetic setting when a change of variables is used (see \eqref{bdry} of Lemma~\ref{lem-CV}).  With this convention, one can easily verify that $F^{-1}_* \eps $ and $F^{-1}_* \mu$ are negative symmetric matrices
since $\det \nabla F(x) < 0$. This clarifies the point that one uses NIMs to construct the complementary medium for the cloaked object. 

Concerning the second part, the medium in $B_{r_1}$ is given by 
\begin{equation}\label{choice-Br1}
\big((r_3^2/r_2^2 ) I, (r_3^2/r_2^2) I \big). 
\end{equation}
%The reason for this choice in $B_{r_1}$  is to fill the space ``disappearing" due to the complementary effect by the homogeneous medium. 
As seen later, mathematically, it is required  to have \eqref{important-identity} in the proof of Theorem~\ref{thm1}. 

\begin{figure}
\centering
\begin{tikzpicture}[scale=1.25]

\draw (0,0) [dashed] circle(2.5);
\draw[fill=black!20] (0,0) circle(1.75);
\draw[fill=red!80] (0,0) circle(1.25);
\draw[fill=blue!80] (0,0) circle(0.75);

\draw[->] (0,0) -- ({0.75*cos(-60)},{0.75*sin(-60)});
\draw ({0.375*cos(-60)},{0.375*sin(-60)}) node[below left]{$r_1$};

\draw[->] (0,0) -- ({1.25*cos(-30)},{1.25*sin(-30)});
\draw ({0.9*cos(-30)},{0.9*sin(-30)}) node[below]{$r_2$};

\draw[->] (0,0) -- (2.5,0);
\draw (2.,0) node[below]{$r_3$};

\draw (0,-2.1) node{$(I,I)$};

\draw [->] (-3.6,2.0) -- ({cos(135)},{sin(135)});
\draw (-3.6,2.0) node[above,align=center]{{\large Complementary layer} \\ $(F^{-1}_* \tilde \varepsilon_O,F^{-1}_* \tilde \mu_O) $};

\draw [->] (4,2.0) -- ({0.375*cos(45)},{0.375*sin(45)});
\draw (4,2.0) node[above,align=center]{{\large Filling space part} \\ $\big((r_3^2/r_2^2)I,(r_3^2/r_2^2)I \big)$};

\draw [->] (-3,-2.0) -- ({1.5*cos(225)},{1.5*sin(225)});
\draw (-3,-2.0) node[below,align=center]{{\large Cloaked object} \\$ (\varepsilon_O, \mu_O) $ };

\end{tikzpicture}
\caption{Cloaking scheme for an object $(\eps_O, \mu_O)$ in $B_{2 r_2} \setminus B_{r_2}$. Two parts are used:  the complementary one in $B_{r_2} \setminus B_{r_1}$ (the red region) which  is the complementary medium of the medium $(\tilde \eps_O, \tilde \mu_O)$ in $B_{r_3} \setminus B_{r_2}$ and the filling space part in $B_{r_1}$ (the blue region)}
\label{fig1}
\end{figure}
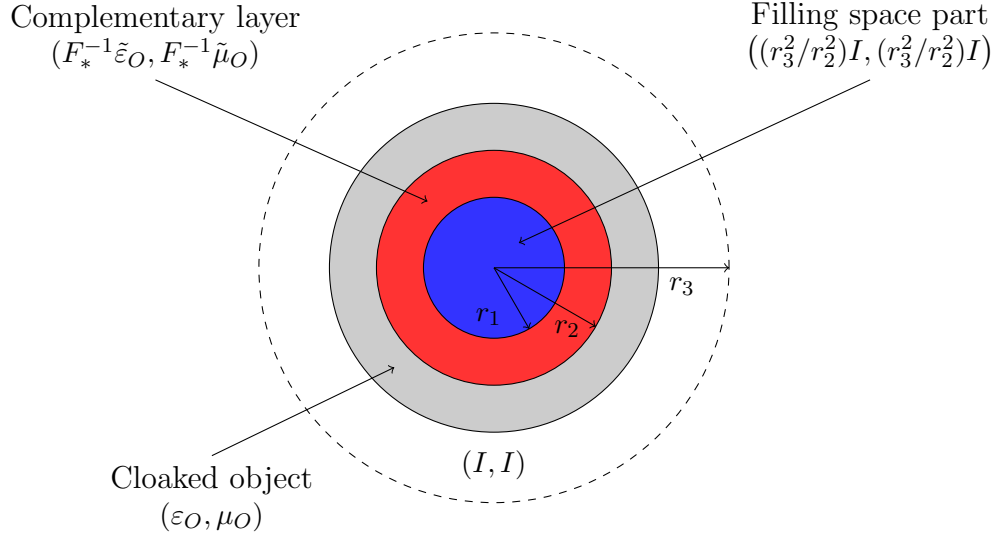

\medskip
Taking into account the loss, the medium in the whole space $\mR^3$ is thus  characterized by $(\eps_\delta,  \mu_\delta)$ 
 defined as follows (see Figure~\ref{fig1} for the case $\delta =0$)
\begin{equation}\label{def-eps-mu}
(\eps_\delta, \mu_\delta) = \left\{ \begin{array}{cl} 
\big( \teps_O , \tmu_O \big)& \mbox{ in } B_{r_3} \setminus B_{r_2}, \\[6pt]
\big( F^{-1}_* \teps_O +  i  \delta I, F^{-1}_* \tmu_O  + i\delta I  \big) & \mbox{ in } B_{r_2} \setminus B_{r_1}, \\[6pt]
\big((r_3^2/r_2^2)  I, (r_3^2/r_2^2)  I \big) & \mbox{ in } B_{r_1},\\[6pt]
 \big(I,  I \big) & \mbox{ in } \mR^3 \setminus B_{r_3}.  
\end{array} \right. 
\end{equation}
Physically, $\eps_\delta$ and $\mu_\delta$ are the permittivity and permeability of the medium, $k$ denotes the frequency, and the imaginary parts of $\eps_\delta$ and $\mu_\delta$ in $B_{r_2} \setminus B_{r_1}$ describe the dissipative property of this (negative index) region. 
% given by a parameter $\delta$. 

%, and $s_0j$ is the current.  The region of negative index material ($\Omega_2 \setminus \Omega_1$) is the region where the real part of $s_\delta$ is negative and   $\delta$ presents the loss of that region.   
\bigskip

Given (a current) $j \in \big[L^2(\mR^3) \big]^3$ with compact support,  let  $(E_\delta, H_\delta), (E, H) \in [H_{\loc}(\curl, \mR^3)]^2$ be respectively the unique outgoing solutions to the Maxwell systems
\begin{equation}\label{eq-EHdelta}
\left\{ \begin{array}{llll}
\nabla \times E_\delta &= & i k  \mu_\delta H_\delta & \mbox{ in } \mR^3\\[6pt]
\nabla \times H_\delta & =&  - i k  \eps_\delta E_\delta + j & \mbox{ in }  \mR^3, 
\end{array} \right.
\end{equation}
and 
\begin{equation}\label{eq-EH}
\left\{ \begin{array}{llll}
\nabla \times E &= & i k   H & \mbox{ in } \mR^3\\[6pt]
\nabla \times H & =&  - i k  E + j & \mbox{ in }  \mR^3. 
\end{array} \right.
\end{equation}

For  an open subset $\Omega$ of $\mR^3$,  the following standard notations  are used: 
\begin{equation*}
H(\curl, \Omega) : = \Big\{ u \in [L^2(\Omega)]^3; \; \nabla \times u  \in [L^2(\Omega)]^3 \Big\}, 
\end{equation*}
\begin{equation*}
\| u\|_{H(\curl, \Omega)} : = \| u\|_{L^2(\Omega)} + \| \nabla \times u \|_{L^2(\Omega)},
\end{equation*}
and 
\begin{equation*}
H_{\loc}(\curl, \Omega) : = \Big\{ u \in [L_{\loc}^2(\Omega)]^3; \; \nabla \times u  \in [L^2_{\loc}(\Omega)]^3 \Big\}. 
\end{equation*}
Recall that a solution $(\cE, \cH) \in [H_{\loc}(\curl, \mR^3 \setminus B_{R})]^2$ (for some $R>0$) of the system 
\begin{equation*}
\left\{ \begin{array}{lll}
\nabla \times \cE &= i k \cH & \mbox{ in } \mR^3 \setminus B_R, \\[6pt]
\nabla \times \cH & = - i k \cE & \mbox{ in }  \mR^3 \setminus B_R, 
\end{array} \right.
\end{equation*}
is said to satisfy the outgoing condition (or the Silver-M\"uller radiation condition) if 
\begin{equation}\label{OC}
\cE \times x + r \cH = O(1/r), 
\end{equation}
as $r = |x| \to + \infty$.

%We assume that 
%\begin{equation}\label{smoothness-1}
%(\teps_O, \tmu_O) \mbox{ is Lipschitz in } B_{r_3} \setminus B_{r_2}. 
%\end{equation}
We extend $(\teps_O, \tmu_O)$ by $(I, I)$ in $B_{r_2}$ and still denote this extension by $(\teps_O, \tmu_O)$. We also assume that 
\begin{equation}\label{smoothness-2}
\mbox{$(\teps_O, \tmu_O)$  is $C^2$ in $B_{r_3}$}.
\end{equation}
Condition~\eqref{smoothness-2} is required for the use of the unique continuation principle and three-sphere inequalities for Maxwell equations. 

%used in the proof of Theorem~\ref{thm1} mainly for the convenience of the presentation, and can be easily fulfilled  by adding an appropriate layer next to the boundary of the cloaked object.   

\medskip 
Cloaking effect of scheme \eqref{def-eps-mu} (see Figure~\ref{fig1}) is mathematically confirmed in the following main result of this paper. 

\begin{theorem}\label{thm1} Let  $R_0 > 0$,   $j \in \big[L^{2}(\mR^3) \big]^3$ with $\supp j \subset \subset B_{R_0} \setminus  B_{r_{3}}$ and let  $(E_\delta, H_\delta),\,  (E, H) \in  \big[H_{\loc}(\curl, \mR^3) \big]^2$ be the unique  outgoing solution to \eqref{eq-EHdelta} and \eqref{eq-EH} respectively. 
Given  $0< \gamma < 1/2$, there exists a positive constant  $\ell = \ell(\gamma)> 0$, depending {\bf only} on  the elliptic  constant of $\teps_O$ and $\tmu_O$ in $B_{2 r_2} \setminus B_{r_2}$ and  $\| (\teps_O, \tmu_O) \|_{W^{2, \infty}(B_{4 r_2})}$ such that if $r_3 > \ell r_2$ then   
\begin{equation}\label{key-point}
\| (E_{\delta}, H_\delta) - (E, H) \|_{H(\curl, B_R \setminus B_{r_3})} \le C_R \delta^\gamma \| j\|_{L^2}, 
\end{equation}
for some positive constant $C_R$ independent of $j$ and $\delta$. 
\end{theorem}

For an observer outside $B_{r_3}$, the medium in $B_{r_3}$ looks like the homogeneous one  by  \eqref{eq-EH}: one has cloaking.

The starting point  of the  proof of Theorem~\ref{thm1} is to use reflections (see \eqref{def-1} and \eqref{def-2}) to obtain Cauchy problems.  We then explore the construction of the cloaking device (its complementary property),  use various three-sphere inequalities (Lemmas~\ref{lem-3S1} and \ref{lem-3S}),  and the removing localized singularity technique to deal with the localized resonance. Using reflections is also the starting point in the study of stability of Helmholtz equations with sign changing coefficients in \cite{Ng-WP} (see also \cite{AnneSophieChesnelCiarlet1,CostabelErnst,Ola} for different approaches) and also plays an important role in the study of superlensing applications of hyperbolic metamaterials in \cite{Bonnetier-Ng}.  A numercial algorithm used for NIMs in the spirit  \cite{Ng-WP} is considered  in \cite{Abdulle16}. Various techniques developed to study NIMs were explored  in the context of interior transmission eigenvalues in \cite{MinhHung1}. The study of NIMs in time domain is recently  investigated in \cite{Joly16}.

\medskip 
The paper is organized as follows. The proof of Theorem~\ref{thm1} is given in Section~\ref{sect-thm1} after presenting several useful results  in Section~\ref{sect-pre}. In Section~\ref{sect-further-discussion}, we present a class of cloaking schemes via the concept of reflecting complementary media.

\section{Preliminaries}\label{sect-pre}

In this section, we present several results which are used in the proof of Theorem~\ref{thm1}.  We first recall a known result  on the trace of $H(\curl, D)$ (see \cite{AlonsoValli, BuffaCostabel}). 

\begin{lemma} \label{lem-trace} Let $D$ be a smooth open bounded subset of $\mR^3$ and set $\Gamma = \partial D$. The tangential trace operator 
\begin{equation*}
\begin{array}{cccc}
\gamma_0 &: H(\curl, D)&   \to &  H^{-1/2}(\dive_\Gamma, \Gamma)\\
& u &  \mapsto & u \times \nu
\end{array}
\end{equation*}
is continuous. Moreover, for all $\phi \in H^{-1/2}(\dive_\Gamma, \Gamma)$, there exists $u \in H(\curl, D)$ such that
\begin{equation*}
\gamma_0(u) = \phi \quad \mbox{ and } \quad \| u\|_{H(\curl, D)} \le C \| \phi\|_{H^{-1/2}(\dive_\Gamma, \Gamma)}, 
\end{equation*}
for some positive constant $C$ independent of $\phi$. 
\end{lemma} 

Here 
\begin{equation*}
H^{-1/2}(\dive_\Gamma, \Gamma): = \Big\{ \phi \in [H^{-1/2}(\Gamma)]^3; \; \phi \cdot \nu = 0 \mbox{ and } \dive_\Gamma \phi \in H^{-1/2}(\Gamma) \Big\}
\end{equation*}
\begin{equation*}
\| \phi\|_{H^{-1/2}(\dive_\Gamma, \Gamma)} : = \| \phi\|_{H^{-1/2}(\Gamma)} +  \| \dive_\Gamma \phi\|_{H^{-1/2}(\Gamma)}.
\end{equation*}

The next result implies the well-posedness and a priori estimates of $(E_\delta, H_\delta)$ defined in \eqref{eq-EHdelta}. 

\begin{lemma}\label{lem-stability-Sys}

Let $k > 0$, $0< \delta < 1$,  $R_0 > 0$,  $D \subset B_{R_0}$
 be  a  smooth bounded open subset of $\mR^3$. Let $\eps, \mu$ be two {\bf real}  measurable matrix-valued functions defined in $\mR^3$ such that  $\eps, \mu$ are uniformly elliptic and piecewise $C^1$ in $\mR^3$, and 
\begin{equation}\label{pro-Identity}
\eps = \mu = I \mbox{ in } \mR^3 \setminus B_{R_0}. 
\end{equation}
Set, for $\delta > 0$, 
\begin{equation}\label{def-eDelta}
(\eps_\delta, \mu_\delta) = \left\{\begin{array}{cl} ( - \eps + i \delta I, - \mu + i \delta I) & \mbox{ if }  x \in D,  \\[6pt]
(\eps, \mu) & \mbox{ otherwise}.
\end{array}\right.
\end{equation}
Let  $j \in L^2(\mR^3)$ with $\supp j \subset B_{R_0}$. There exists  a unique outgoing solution $(E_\delta, H_\delta) \in [H_{\loc}(\curl, \mR^3)]^2$  to the Maxwell system
\begin{equation}\label{eq-EHDelta}
\left\{ \begin{array}{llll}
\nabla \times \cE_\delta &= & i k \mu_\delta \cH_\delta & \mbox{ in } \mR^3,\\[6pt]
\nabla \times \cH_\delta & =&  - i k \eps_\delta \cE_\delta + j & \mbox{ in }  \mR^3. 
\end{array} \right.
\end{equation}
Moreover, 
\begin{equation}\label{stability2}
\|(\cE_\delta, \cH_\delta) \|_{H(\curl, B_{R})}^2 \le C_R \left(  \frac{1}{\delta} \| j\|_{L^2} \|(\cE_\delta, \cH_\delta) \|_{L^2(\supp j)} + \|j \|_{L^2}^2 \right). 
\end{equation}
Here $C_R$ denotes a positive constant depending on $R$, $R_0$, $\eps$, $\mu$ but independent of $j$ and $\delta$. Consequently, we have 
\begin{equation}\label{stability1-11}
\|(\cE_\delta, \cH_\delta) \|_{H(\curl, B_{R})} \le \frac{C_{R}}{\delta} \| j \|_{L^2}. 
\end{equation}
\end{lemma}

%Lemma~\ref{lem-stability-Sys} is a direct consequence of \cite[Lemma 6]{Ng-Superlensing-Maxwell} by taking $f = 0$, $g = j$, $\Omega_{2} = B_{r_2}$, $\Omega_1  = B_{r_1}$, and $D = B_{R_0} \setminus B_{r_3}$. 

\noindent{\bf Proof.} The existence of $(\cE_\delta, \cH_\delta)$ can be derived from the uniqueness of $(\cE_\delta, \cH_\delta)$ as usual. The uniqueness of $(\cE_\delta, \cH_\delta)$ can be deduced from the estimates of $(\cE_\delta, \cH_\delta)$. Estimate~\eqref{stability1-11} is a direct consequence of \eqref{stability2}. We hence only give the proof of  \eqref{stability2}. We have, by \eqref{eq-EHDelta},  
\begin{equation*}
\nabla \times ( \mu_\delta^{-1} \nabla \times \cE_\delta) - k^2 \eps_\delta \cE_\delta =   i k j \mbox{ in } \mR^3. 
\end{equation*}
Set 
\begin{equation*}
M_\delta = \frac{1}{\delta} \| j \|_{L^2} \| (\cE_\delta, \cH_\delta) \|_{L^2(\supp j)} + \|j\|_{L^2}^2. 
\end{equation*}
Multiplying the equation by $\bar \cE_\delta$, integrating in $B_R$, and using the fact that $\supp j \subset B_{R_0}$, we have, for $R > R_0$, 
\begin{equation*}
\int_{B_R} \langle  \mu_\delta^{-1} \nabla \times \cE_\delta, \nabla \times \cE_\delta \rangle  - \int_{\partial B_{R}}  \big\langle (\mu_\delta^{-1} \nabla \times \cE_\delta) \times \nu, \cE_\delta \big\rangle  - k^2 \int_{B_R} \langle \eps_\delta \cE_\delta, \cE_\delta \rangle =   \int_{B_R}  \langle i k j, \cE_\delta \rangle.
\end{equation*}
Since  $\mu_{\delta} = I$ and so $\nabla \times  \cE_\delta = i k \cH_\delta $ in $\mR^3 \setminus B_{R_0}$, we derive that, for $R> R_0$,  
\begin{equation*}
\int_{B_R} \langle  \mu_\delta^{-1} \nabla \times \cE_\delta, \nabla \times \cE_\delta \rangle  + \int_{\partial B_{R}}  \langle i k \cH_\delta , \cE_\delta \times \nu \rangle  - k^2 \int_{B_R} \langle \eps_\delta \cE_\delta, \cE_\delta \rangle  =    \int_{B_R}  \langle i k j, \cE_\delta \rangle.
\end{equation*}
Letting $R \to + \infty$, using the outgoing condition ($\cE_\delta(x) \times \nu(x) = - \cH_\delta(x) + O(1/R^2)$ for $x \in \partial B_R$),  and considering the imaginary part,  we obtain 
\begin{equation}\label{inE}
 \|\cE_\delta \|_{H(\curl, D)}^2 \le C M_\delta.    
\end{equation}
This  implies, by Lemma~\ref{lem-trace},  with the notation $\Gamma = \partial D$, 
\begin{equation}\label{bdryE}
\| \cE_\delta \times \nu \|_{H^{-1/2}(\dive_\Gamma, \Gamma)}^2   \le C M_\delta. 
\end{equation}
Using the equations of $(\cE_\delta, \cH_\delta)$ in $D$, we derive from \eqref{inE} that 
\begin{equation}\label{inH}
 \|\cH_\delta \|_{H(\curl, D)}^2 \le C M_\delta; 
\end{equation}
which yields, by Lemma~\ref{lem-trace} again,  
\begin{equation}\label{bdryH}
\| \cH_\delta \times \nu \|_{H^{-1/2}(\dive_\Gamma,\Gamma)}^2  \le C M_\delta. 
\end{equation}
Let $D_{1}^c$ be the unbounded connected component of $\mR^3 \setminus \bar D$ and let  $D_2^c$ be the complement of $D_{1}^c$ in $\mR^3 \setminus \bar D$, i.e., $D_2^c = (\mR^3 \setminus \bar D) \setminus D_{1}^c$ \footnote{We will apply Lemma~\ref{lem-stability-Sys} with $D = B_{r_2} \setminus B_{r_1}$; in this case $D_1^c = \mR^3 \setminus \bar B_{r_2}$ and $D_2^c = B_{r_1}$.}. We have
\begin{equation*}
\left\{ \begin{array}{llll}
\nabla \times \cE_\delta &= & i k \mu \cH_\delta & \mbox{ in } D_1^c,\\[6pt]
\nabla \times \cH_\delta & =&  - i k \eps \cE_\delta + j & \mbox{ in }  D_1^c. 
\end{array} \right.
\end{equation*}
It follows that, see  e.g., \cite[Lemma 5]{Ng-Superlensing-Maxwell},  
\begin{equation}\label{inEHo}
\| (\cE_\delta, \cH_\delta)\|_{H(\curl, B_R \cap D_1^c)}^2 \le C_R \big( \|j \|_{L^2} + \|(\cE_\delta, \cH_\delta)\|_{H^{-1/2}(\dive_{\partial D_1^c}, \partial D_1^c)}  \big),
\end{equation}
We deduce from \eqref{bdryE}  that 
\begin{equation}\label{inEHo}
\| (\cE_\delta, \cH_\delta)\|_{H(\curl, B_R \cap D_1^c)}^2 \le C_R M_\delta,
\end{equation}
and, by  Lemma~\ref{lem-stability-inside} below, we derive from \eqref{bdryE} and \eqref{bdryH} that 
\begin{equation}\label{inEHi}
\| (\cE_\delta, \cH_\delta)\|_{H(\curl,  D_2^c)}^2 \le C M_\delta. 
\end{equation}
A combination of \eqref{inE}, \eqref{inH}, \eqref{inEHo}, and \eqref{inEHi} yields 
\begin{equation}\label{hahahaha}
\| (\cE_\delta, \cH_\delta)\|_{H(\curl, B_R)} \le C_R M_\delta; 
\end{equation}
which is  \eqref{stability2}.   \proofend

\medskip

In the proof of Lemma~\ref{lem-stability-Sys} we use the following result (see \cite[Lemma 3]{Ng-Superlensing-Maxwell}) whose proof follows directly from  the unique continuation principle for the Maxwell equations (see \cite{BallCapdeboscq, Tu}) via a contradiction argument. 

\begin{lemma}\label{lem-stability-inside}
 Let $k>0$, $D $ be a smooth bounded open subset of $\mR^3$, $f, g \in [L^2(D)]^3$, and $h_1, h_2 \in H^{-1/2}(\dive_{\partial D}, \partial D)$, and let $\eps$ and $\mu$ be two piecewise $C^1$, symmetric uniformly elliptic matrix-valued functions defined in $D$. Assume that  $(\cE, \cH) \in [H( \curl, D)]^2  $ is a  solution to
\begin{equation*}
\left\{ \begin{array}{cll}
\nabla \times \cE = i k  \mu \cH + f& \mbox{ in } D, \\[6pt]
\nabla \times \cH  = - i k \eps \cE + g & \mbox{ in }  D, \\[6pt]
\cH \times \nu = h_1; \; \cE \times \nu = h_2 & \mbox{ on } \partial D.  
\end{array} \right.
\end{equation*}
Then
\begin{equation}\label{stability1-1}
\|(\cE, \cH) \|_{H(\curl, D)} \le C \Big(\|(f, g) \|_{L^2(D)} +  \| (h_1, h_2)\|_{H^{-1/2}(\dive_\Gamma, \partial D)}\Big), 
\end{equation}
for some positive constant $C$ depending on $D$, $\eps$, $\mu$, and $k$ but independent of $f$, $g$, $h_1$, and $h_2$. 
\end{lemma}

We next state a  result  on the existence, uniqueness, and stability results on Maxwell equations with elliptic coefficients (see \cite[Lemma 4]{Ng-Superlensing-Maxwell}).

\begin{lemma} \label{lem-stability} Let $k>0$, $D$ be a smooth  bounded open subset of $\mR^3$, $f, g \in [L^2(\mR^3)]^3$, $h_1, h_2 \in  H^{-1/2}(\dive_{\partial D}, \partial D)$. Assume that $\bar D, \supp f, \supp g \subset B_{R_0}$ for some $R_0> 0$. Let $\eps, \mu$ 
be two symmetric uniformly elliptic matrix-valued functions defined in $\mR^3$ such that $(\eps, \mu) = (I, I)$ in $\mR^3 \setminus B_{R_0}$ and $(\eps, \mu)$ is piecewise $C^1$. There exists  a unique solution $(\cE, \cH) \in \big[\bigcap_{R>0} H(\curl, B_R \setminus \partial D) \big]^2$ of the system 
\begin{equation}\label{haha}
\left\{ \begin{array}{cl}
\nabla \times \cE = i k  \mu \cH + f & \mbox{ in } \mR^3 \setminus \partial D, \\[6pt]
\nabla \times \cH  = - i k \eps \cE + g & \mbox{ in }  \mR^3 \setminus \partial D, \\[6pt]
[\cH \times \nu] = h_1; \; [\cE \times \nu] = h_2 & \mbox{ on } \partial D.  
\end{array} \right.
\end{equation}
Moreover, 
\begin{equation}\label{stability1}
\|(\cE, \cH) \|_{  H(\curl, B_{R} \setminus \partial D)} \le C_{R} \Big(\| (f, g) \|_{L^2} + \| (h_1, h_2) \|_{H^{-1/2}(\dive_\Gamma, \partial D)} \Big), 
\end{equation}
for some positive constant $C_R$ depending on $R$, $R_0$, $D$, $\eps$, $\mu$, and $k$,  but independent of $f, g, h_1$, and $h_2$. 
\end{lemma}

Hereafter $[ \cdot ]$ denotes the jump across the boundary $\partial D$, i.e., $[u] = u|_{\mathrm{ext}} - u|_{\mathrm{int}}$ for an appropriate function $u$. 

\begin{remark}  \fontfamily{m} \selectfont
%The well-posedness of \eqref{haha} is known for $h_1= h_2 = 0$ and $f, g \in H(\dive, \mR^3)$ where $H(\dive, \mR^3): = \{u \in [L^2(\mR^3)]^3; \; \dive u \in L^2(\mR^3)\}$ and $\| u\|_{H(\dive, \mR^3)} : = \| u \|_{L^2(\mR^3)} + \| \dive u\|_{L^2(\mR^3)}$. In this case, $\| (f, g) \|_{L^2} $  is replaced by $\| (f, g) \|_{H(\dive, \mR^3)}$ in \eqref{stability1}. 
In the proof of Lemma~\ref{lem-stability}, one uses a new compactness criterion for the Maxwell equations in \cite{Ng-Superlensing-Maxwell} with roots in \cite{HaddarJolyNguyen2}.  
\end{remark}

%\begin{remark}  \fontfamily{m} \selectfont  System \eqref{interior-Lem} is known for a smooth pair $(E, H)$.   %(see, e.g., \cite{PendrySchurigSmith}). 
%Statement \eqref{bdry}  might be known; however, we cannot find a reference for it.  For the convenience of the reader, we give  the details of the proof in the appendix for the  form stated here. 
%\end{remark}

We next present a known result  which reveals a connection between Maxwell equations with weakly coupled elliptic systems. 

\begin{lemma}\label{lem-elliptic} Let $D$ be an open subset of $\mR^3$, $\eps, \mu$ be two matrix-valued functions defined in $D$,  and let $(\cE, \cH) \in \big[H^1(D) \big]^2$ be  a solution of the system
\begin{equation}\label{sys-elliptic}
\left\{\begin{array}{cl}
\nabla \times \cE = i  k \mu \cH & \mbox{ in } D, \\[6pt]
\nabla \times \cH = - i k \eps \cE & \mbox{ in } D. 
\end{array}\right.
\end{equation}
Then, for $1 \le a \le 3$,  
\begin{equation}\label{eq-H}
\dive (\mu \nabla \cH_a)  + \dive (\partial_a \mu \cH - i k \mu \epsilon^a \eps \cE) = 0 \mbox{ in } D, 
\end{equation}
\begin{equation}\label{eq-E}
\dive(\eps \nabla \cE_a) + \dive (\partial_a \eps \cE + i k \eps \epsilon^a \mu \cH) = 0 \mbox{ in } D. 
\end{equation}
Here $\epsilon^a_{bc}$ ($1 \le a, b, c \le 1$) denotes the usual Levi Civita permutation, i.e., 
\begin{equation}\label{Levi-Civita}
\epsilon^a_{bc} = \left\{\begin{array}{cl} \mbox{sign } (abc) & \mbox{ if  $abc$  is a permuation}, \\[6pt]
0 & \mbox{otherwise}. 
\end{array}\right.
\end{equation}
\end{lemma}

\noindent{\bf Proof.}  The proof is quite simple as follows. Using the fact, for $1 \le a \le 3$,  
$$
\partial_a \cH = \nabla \cH_a + \epsilon^a (\nabla \times \cH) \quad  \mbox{ and } \quad \partial_a \cE = \nabla \cE_a + \epsilon^a (\nabla \times  \cE),  
$$
we derive from \eqref{sys-elliptic} that, for $1 \le a \le 3$,  
\begin{equation}\label{elliptic-part1}
\partial_a \cH  = \nabla \cH_a - i k \epsilon^a \eps \cE \quad \mbox{ and }  \quad  \partial_a \cE = \nabla \cE_a + i k \epsilon^a \mu \cH \mbox{ in } D. 
\end{equation}
Since
\begin{equation*}
 \dive (\mu \cH) = 0 \mbox{ in } D, 
\end{equation*}
it follows that, for $1 \le a \le 3$,  
\begin{equation*}
0 =  \partial_a \dive (\mu \cH) = \dive (\mu \partial_a \cH) + \dive(\partial_a \mu \cH) \mbox{ in } D.
\end{equation*}
This implies, by the first identity of \eqref{elliptic-part1},  
\begin{equation*}
\dive (\mu \nabla \cH_a)  + \dive (\partial_a \mu \cH - i k \mu \epsilon^a \eps \cE) = 0 \mbox{ in } D; 
\end{equation*}
which is \eqref{eq-H}.  Similarly, we obtain \eqref{eq-E}.  \proofend

\medskip 
Hadamard proved the following three-circle inequality: Assume that $\Delta v = 0$ in $B_{R^*} \setminus B_{R_*} \subset \mR^2$ and $0< R_* < R_1 < R_2 < R_3 < R^*$. Then 
$$
\| v\|_{L^\infty(\partial B_{R_2})} \le \| v\|_{L^\infty(\partial B_{R_1})}^\alpha \| v\|_{L^\infty(\partial B_{R_3})}^{1 - \alpha}, 
$$
with $\alpha =  \ln (R_3/ R_2) / \ln(R_3/ R_1)$. Here is its variant  which is useful in the proof of Theorem~\ref{thm1}. 

\begin{lemma} \label{lem-3S1} Let $d=2, \, 3$, $k, \, R_*, \, R^*>0$,   and let $v \in H^1(B_{R_*} \setminus B_{R^*})$ be a solution to the equation $\Delta v + k^2 v = 0$ in $B_{R_3} \setminus B_{R_1} \subset \mR^d$. We have, for $R_* \le R_1 < R_2 < R_3 \le R^*$, 
\begin{equation}\label{ineq-3S}
\| v \|_{\bH(\partial B_{R_2})} \le C \| v\|_{\bH(B_{R_1})}^{\alpha} \| v \|_{\bH(B_{R_3})}^{1 - \alpha}, 
\end{equation}
where $\alpha =  \ln (R_3/ R_2) / \ln(R_3/ R_1)$ and $C$ is a positive constant depending only on $k$, $R_*$, and $R^*$. Here 
$$
\| v\|_{\bH(\partial B_r)}: = \|v \|_{H^{1/2}(\partial B_r)} + \|\partial_r v \|_{H^{-1/2}(\partial B_r)}. 
$$
\end{lemma}

%
%\begin{remark} \fontfamily{m} \selectfont  Inequality \eqref{ineq-3S} is a variant of the well-known Hadamard's three-sphere inequality in which $\| v\|_{\bH(\partial B_r)}$ is replaced by $\|v \|_{L^\infty(\partial B_r)}$, and  $d=2$ and $k=0$;   the exponent $\alpha$ in Lemma~\ref{lem-3S1} is the same as in Hadamard's three-sphere inequality. 
%\end{remark}

%\begin{remark} \fontfamily{m} \selectfont  The case $k=0$ is well-known and first noted by Hadamard in two dimensions (see \cite{MinhLoc2} for a recent discussion on three spheres inequalities for second order elliptic equations and their applications for cloaking using complementary media). 
%\end{remark}

\medskip 
Before giving the proof of Lemma~\ref{lem-3S1}, we recall some properties of the spherical Bessel and Neumann functions and the Bessel and Neumann functions of large order. We first introduce,  for $n \ge 1$, 
\begin{equation}\label{def-jn}
\hat j_n(t) =1 \cdot 3 \cdots (2n + 1) j_n(t) \quad \mbox{ and } \quad  \hat y_n = -  \frac{y_n(t)}{1 \cdot 3 \cdots (2n-1)} ,  
\end{equation}
and  for $n \ge 0$, 
\begin{equation}\label{def-Jn}
\hat J_n(r) = 2^n n! J_n(r) \quad \mbox{ and } \quad \hat Y_n(r) = \frac{\pi i}{2^{n} (n-1)!} Y_n(r), 
\end{equation}
where $j_n$ and $y_n$   are the spherical Bessel and Neumann functions, and 
 $J_n$ and $Y_n$ are the Bessel and Neumann functions  of order $n$ respectively. Then, see, e.g.,  \cite[(2.37), (2.38), (3.57), and  (3.58)]{ColtonKressInverse}),   as $n \to + \infty$,  
\begin{equation}\label{jy-n}
\hat j_n(r) = r^n \big[1 + O(1/n) \big] \quad \hat y_n(r) = r^{-n-1} \big[1 + O(1/n) \big], 
\end{equation}
\begin{equation}\label{JY-n}
\hat J_n (t)  = t^{n}\big[1 + O(1/n) \big], \quad \mbox{ and } \quad  \hat Y_n (t)  = t^{-n} \big[1 + O(1/n) \big]. 
\end{equation}
One also has, see, e.g.,   \cite[(2.36) and (3.56)]{ColtonKressInverse},
\begin{equation}\label{W3}
  j_n(r)  y_n'(r) -   j_n'(r) y_n(r) = \frac{1}{ r^2}
\end{equation}
and 
\begin{equation}\label{W2}
  J_n(r)  Y_n'(r) -   J_n'(r) Y_n(r) = \frac{2}{ \pi r}. 
\end{equation}

%Using \eqref{bh1} and \eqref{jy-n}, we can now implement the analysis in the quasistatic regime developed in \cite{Ng-CALR} to the finite frequency regime in this section. 
%Similar ideas were used in \cite{Ng-Superlensing, MinhLoc2}. 

We are ready to give 

\medskip
\noindent{\bf Proof of Lemma~\ref{lem-3S1}.} By rescaling, one can assume that $k =1$.  We consider the case $d=2$ and $d=3$ separately. 

\medskip 
\noindent {\bf Case 1: d=3.} Since $\Delta v + v = 0$ in $B_{R_3} \setminus B_{R_1}$, $v$ can be represented in the form 
\begin{equation*}
v=  \sum_{n=1}^\infty \sum_{m= -n }^n \Big( a^n_m \hat j_{n}(|x|)  + b^n_m \hat y_{n}(|x|)  \Big) Y^n_m(\hat x ) \quad \mbox{ in } B_{R_3} \setminus B_{R_1}, 
\end{equation*}
for $ a^n_m \in \mC$ and $\hat x = x/ |x|$ where $Y^n_m$ is the spherical harmonic function of degree $n$ and of order $m$. In what follows in this proof, $C$ denotes a positive constant depending only on  $R_*$ and $R^*$ and can change from one place to another and $a \sim b$ means that $a \le Cb$ and $b \le Ca$.  Using the fact $(Y^{n}_m)$ is an orthonormal basis of $L^2(\partial B_1)$ and 
$$
\Delta_{\partial B_1} Y^n_m + n (n+1) Y^n_m = 0 \mbox{ on } \partial B_1,
$$ 
we derive that, for $R_1 \le r \le R_3$,
\begin{equation}\label{norm-equivalent-3}
\| v\|_{\bH(\partial B_r)}^2 \sim \sum_{n =1}^\infty \sum_{m = - n}^n \big(  n   |c^n_m(r)|^2  + n^{-1} |d^n_m(r)|^2 \big),
\end{equation} 
where 
\begin{equation}\label{cn-dn}
c^n_m(r) =  a^n_m \hat j_{n}(r)  + b^n_m \hat y_{n}(r) \quad \mbox{ and } \quad d^n_m(r) = a^n_m \hat j_{n}'(r)  + b^n_m \hat y_{n} ' (r). 
\end{equation}
From \eqref{norm-equivalent-3} and \eqref{cn-dn}, we have
\begin{equation*}
\| v\|_{\bH(\partial B_r)}^2 \le C \sum_{n =1}^\infty \sum_{m = - n}^n  \Big(  |a^n_m|^2 \big(n  |\hat j_{n}(r)|^2  + n^{-1}  |\hat j_{n}'(r)|^2 \big) +   |b^n_m|^2 \big(n  |\hat y_{n}(r)|^2  + n^{-1}  |\hat y_{n}'(r)|^2 \big)  \Big); 
\end{equation*}
which yields, by \eqref{jy-n},  
\begin{equation}\label{norm-equivalent-3-1}
\| v\|_{\bH(\partial B_r)}^2 \le  C \sum_{n =1}^\infty \sum_{m = - n}^n \big( n r^{2n} |a^n_m|^2  + n r^{-2n} |b^n_m|^2 \big). 
\end{equation} 
From \eqref{cn-dn}, we have 
\begin{equation}\label{an-bn}
a^n_m = \frac{c^n_m(r) \hat y_n'(r) - d^n_m(r) \hat y_n(r) }{  \hat j_n(r) \hat y_n'(r) - \hat j_n'(r) \hat y_n(r)} \quad \mbox{ and } \quad b^n_m = \frac{c^n_m(r) \hat j_n'(r) - d^n_m(r) \hat j_n(r) }{  \hat y_n(r) \hat j_n'(r) - \hat y_n'(r) \hat j_n(r)}. 
\end{equation}
From \eqref{W3}, we obtain, for some $c_n \neq 0$, 
\begin{equation}\label{W3-1}
\hat y_n(r) \hat j_n'(r) - \hat y_n'(r) \hat j_n(r) = \frac{c_n}{t^2}. 
\end{equation}
Combining  \eqref{jy-n}, \eqref{an-bn}, and \eqref{W3-1} yields  
\begin{equation}\label{an-bn-cn-dn}
|a^n_m| \le C \big(  |c^n_m| r^{-n} + n^{-1} |d^n_m| r^{-n} \big) \quad \mbox{ and } \quad |b^n_m| \le C \big( |c^n_m| r^{n} + n^{-1} |d^n_m| r^n \big). 
\end{equation}
We derive from  \eqref{norm-equivalent-3} and \eqref{an-bn-cn-dn} that 
\begin{equation}\label{norm-equivalent-3-2}
 \sum_{n =1}^\infty \sum_{m = - n}^n \big( n r^{2n} |a^n_m|^2  + n r^{-2n} |b^n_m|^2 \big) \le C  \| v\|_{\bH(\partial B_r)}^2. 
\end{equation}
A combination of \eqref{norm-equivalent-3-1} and \eqref{norm-equivalent-3-2} yields 
\begin{equation}\label{norm-equivalent-3-3}
\| v\|_{\bH(\partial B_r)}^2 \sim   \sum_{n =1}^\infty \sum_{m = - n}^n \big( n r^{2n} |a^n_m|^2  + n r^{-2n} |b^n_m|^2 \big). 
\end{equation}
Inequality \eqref{ineq-3S} is now a consequence of \eqref{norm-equivalent-3-3}  after applying H\"older's inequality and noting that $R_2 = R_1^\alpha R_3^{1 - \alpha}$.

\medskip 
\noindent {\bf Case 2: d=2.} Since $\Delta v + v = 0$ in $B_{R^*} \setminus B_{R_*}$, one can represent $v$ of the form
\begin{equation*}
v= a_0 + b_0 \ln r +   \sum_{n=1}^\infty \sum_{\pm} \big( a_{n, \pm} \hat J_{n}(|x|) + b_{n, \pm} \hat Y_{n}(|x|) \big)e^{\pm i n \theta}  \mbox{ in } B_{R^*} \setminus B_{R_*}, 
\end{equation*} 
Using \eqref{JY-n} and \eqref{W2}, as in the previous case, one can prove that 
\begin{equation}
\| v\|_{\bH(\partial B_r)}^2 \sim |a_0|^2 + |b_0|^2 + \sum_{n \ge 1} \sum_{\pm} \big( n r^{2n} |a_{n, \pm}|^2 + n^{-1} r^{-2n} |b_{n, \pm}|^2\big)
\end{equation}
Inequality \eqref{ineq-3S} is now a consequence of \eqref{norm-equivalent-3-3}  after applying H\"older's inequality and noting that $R_2 = R_1^\alpha R_3^{1 - \alpha}$.  \proofend

\medskip 
We next state a three-sphere inequality for an ``elliptic system". 

\begin{lemma}\label{lem-3S}
Let $m, n \in \mN$ $(m \ge 2, \; n \ge 1)$, $0 < R_* <  R_1 < R_2 < R_3 < R^*$, $c_1, c_2 > 0$  and let $M^1, \cdots, M^n$ be such that $M^k$ is an $(m \times m)$ matrix defined in $B_{R^*} \subset \mR^m$  for $1 \le k \le n$ \footnote{In this lemma, $B_r$ denotes the ball centered at the origin with radius $r$ in $\mR^m$.}.  Assume that $M^k$ is Lipschitz and  uniformly elliptic in $B_{R^*}$  for $1 \le k \le n$ and   $V \in \big[H^1(B_{R_3} \setminus \bar B_{R_1}) \big]^n$ satisfies 
\begin{equation}\label{eq-V}
|\dive (M^k \nabla V_k)| \le  c_{1} |\nabla V| + c_2 |V|  \mbox{ a.e. in } B_{R_3} \setminus \bar B_{R_1} \mbox{ for } 1 \le k \le n. 
\end{equation}
There exists a constant $q \ge 1$, depending only on $m, \, n$, and the elliptic and the Lipschitz constants of $M^k$ for $1 \le k \le n$ such that, for any $\lambda_0 > 1$ and $R_2 \in (\lambda_ 0 R_1, R_3/\lambda_0)$, we have 
 \begin{equation}\label{def-alpha}
\| V\|_{\bH(\partial B_{R_2})} \le C\| V\|_{\bH(\partial B_{R_1})} ^{\alpha} \| V\|_{\bH(\partial B_{R_3})}^{1-\alpha}
\quad \mbox{ where } \quad 
	\alpha := \frac{R_2^{-q} - R_3^{-q}}{R_1^{-q} - R_{3}^{-q}}, 
\end{equation}
and 
\begin{equation}\label{H-norm}
\| V_k \|_{\bH(\pB_r)} = \| V_k\|_{H^{1/2}(\pB_r)} + \| M_k \nabla V_k \cdot e_r \|_{H^{-1/2}(\pB_r)}, \quad \|V \|_{{\bH(\pB_r)}} = \sum_{k=1}^n \| V_k \|_{\bH(\pB_r)}. 
\end{equation} 
Here $C$ is a positive constant  depends on the elliptic and the Lipschitz constants of $M^k$ $(1 \le k \le n)$, $c_1,$ $c_2,$  $R_*, R^*,  m, n$, and $\lambda_0$ but independent of $v$. 
\end{lemma}

In inequality \eqref{def-alpha}, the constant $q$ does {\bf not} depend on $c_1, \, c_2, R_*, R^*$ but the constant $C$  {\bf does}. No upper bound on $R^*$ is imposed as often required in a three-sphere inequality for the Helmholtz equation (see e.g.,  \cite[Theorem 4.1]{AlessandriniRondi}). Nevertheless, both information of $V$ and $M \nabla V \cdot e_r$ are used \eqref{H-norm}; this is the key point to ensure that \eqref{def-alpha} holds without imposing any condition on $R^*$. Lemma~\ref{lem-3S} is proved in \cite[Theorem 2]{MinhLoc2} for the case $n=1$. The proof for the case $n \ge 1$ follows similar and is omitted. 

\medskip 
We finally state a change of variables formula,  which can be found in \cite[Lemma 7]{Ng-Superlensing-Maxwell}. 

\begin{lemma} \label{lem-CV}Let $D, D'$ be two bounded connected open subsets of $\mR^3$ and ${\cal T}: D \to D'$ be  bijective such that 
${\cal T} \in C^1(\bar D)$ and ${\cal T}^{-1} \in C^1(\bar D')$. Assume that $\eps, \mu \in [L^\infty(D)]^{3 \times 3}$,  $j \in [L^2(D)]^3$ and $(E, H) \in [H(\curl, D)]^2$ is a solution to 
\begin{equation*}\left\{
\begin{array}{lll}
\nabla \times E  &=  i k \mu H & \mbox{ in } D, \\[6pt]
\nabla \times H &=  - i k \eps  E + j & \mbox{ in } D. 
\end{array}\right. 
\end{equation*}
Define $(E', H')$ in $D'$ as follows 
\begin{equation}\label{TEH}
E'(x') = {\cal T}*E(x'): = \nabla {\cal T}^{-T}(x) E(x) \mbox{  and  } H'(x') = {\cal T}*H(x'): = \nabla {\cal T}^{-T}(x) H(x),  
\end{equation}
with $x' = {\cal T}(x)$.   Then $(E', H')$ is a solution to 
\begin{equation}\label{interior-Lem}
\left\{
\begin{array}{lll}
\nabla' \times E'  &=  i k \mu' H' & \mbox{ in } D', \\[6pt] 
\nabla' \times H' &=  - i k \eps'  E' + j' & \mbox{ in } D', 
\end{array}\right. 
\end{equation}
where 
\begin{equation*}
\eps' = {\cal T}_*\eps, \quad \mu' = {\cal T}_*\mu, \quad j'= T_*j,
\end{equation*}
and 
\begin{equation}\label{def-T**}
{\cal T}_* \eps (x') = \frac{\nabla {\cal T}  (x)  \eps(x) \nabla  {\cal T} ^{T}(x)}{J(x)}, \quad  {\cal T}_* \mu (x') = \frac{\nabla {\cal T}  (x)  \mu(x) \nabla  {\cal T} ^{T}(x)}{J(x)},  \quad \mbox{ and } \quad {\cal T}_*  j (x') = \frac{j(x)}{J(x)},
\end{equation}
with $x ={\cal T}^{-1}(x')$ and $J(x) = \det \nabla  {\cal T}(x)$. 
Assume in addition that  $D$ is of class $C^1$ and ${\bf T} = {\cal T} \big|_{\partial D}: \partial D \to \partial D'$ is a diffeomorphism. We have \footnote{Here $\nu$ and $\nu'$ denote the outward unit normal  vector on $\partial D$ and $\partial D'$.} 
\begin{equation}\label{bdry}
\mbox{if } E \times \nu = g  \mbox{ and } H \times \nu = h \mbox{ on } \partial D \mbox{ then } E' \times \nu' = {\bf T}_*g \mbox{ and } H' \times \nu' = {\bf T}_* h  \mbox{ on } \partial D', 
\end{equation}
where ${\bf T}_*$ is given in \eqref{def-T} below. 
\end{lemma}

%Given a diffeomorphism $T$ from $D$ onto $D'$, the following standard notations are used 
%\begin{equation}\label{def-TT}
%{\cal T}_*a (x') = \frac{\nabla {\cal T}(x) a(x) \nabla {\cal T}^T(x)}{\det \nabla {\cal T}(x)}  \quad \mbox{ and } \quad  {\cal T}_*j (x')= \frac{\nabla {\cal T}(x) j(x)}{\det \nabla {\cal T}(x)},
%\end{equation}
%with  $x' ={\cal  T}(x)$,  for a matrix-valued function $a$ and   a vector-valued function $j$ defined in $D$.  
For a tangential vector field $g$  defined in $\partial D$,  we denote
\begin{equation}\label{def-T}
{\bf T}_*g(x') = \sign \cdot \frac{\nabla_{\partial D} {\bf T}(x) g(x)}{|\det \nabla_{\partial D}  {\bf T}(x)|}  \mbox{ with } x' = {\bf T}(x),  
\end{equation}
where $\sign := \det \nabla {\cal T} (x)/ |\det \nabla {\cal T}(x)|$ for some $x \in D$.

\begin{remark} \fontfamily{m} \selectfont
In the change of variables, the definition of ${\cal T}*$ in \eqref{TEH} is different from ${\cal T}_*$ in \eqref{def-T**}. It is worthy remembering that for electromagnetic fields  \eqref{TEH} is used whereas for sources \eqref{def-T*} is involved. %  In the definition of $T_*$ given in \eqref{def-T**}, $J(x): = \det \nabla {\cal T}(x) $ {\bf not} $|\det \nabla {\cal T}(x)|$ as often used in the acoustic setting. 
\end{remark}

Lemma~\ref{lem-CV} places an important role in the construction of the cloaking device and is the ingredient for the reflecting technique used in the proof of Theorem~\ref{thm1}.

\section{Proof of Theorem~\ref{thm1}}\label{sect-thm1}

%It is known that resonance might occur in a medium which has complementary property. This was established in \cite[Proposition 2]{Ng-WP} for the Helmholtz equation; the situation for the Maxwell equations would hold with similar figure. To deal with resonance, and more precisely, localized resonance, we introduce the removing localized sigularity technique for the Maxwell equations. This technique was introduced in \cite{Ng-Negative-Cloaking, Ng-Superlensing} in the acoustic setting.  The proof of Theorem~\ref{thm1} is also based on the reflecting technique via Lemma~\ref{lem-CV} and a three-sphere inequality via Lemma~\ref{lem-3S}. 
%
%
%\medskip 
%\noindent{\bf Proof of Theorem~\ref{thm1}.} 

Let $(\oE_\delta, \oH_\delta) \in \big[H^1_{\loc}(\curl, \mR^3 \setminus B_{r_2}) \big]^2$ be the reflection of $(E_\delta, H_\delta)$  through $\partial B_{r_2}$ by the Kelvin transform $F$ with respect to $\partial B_{r_2}$, i.e., 
\begin{equation}\label{def-1}
\big(\oE_\delta, \oH_\delta \big)= \big( F*E_\delta, F*H_\delta \big) \mbox{ in } \mR^3 \setminus B_{r_2},  
\end{equation}
where $F*$ is defined by \eqref{TEH}.  Let $\big(\tE_\delta, \tH_\delta \big) \in \big[H(\curl, B_{r_3}) \big]^2$ be the reflection of $(\oE_\delta, \oH_\delta)$  through $\partial B_{r_3}$ by  the Kelvin transform $G: \mR^3 \setminus B_{r_3} \mapsto B_{r_3}$ with respect to $\partial B_{r_3}$., i.e., $G(x) = r_3^2 x / |x|^2 $ and 
\begin{equation}\label{def-2}
\big(\tE_\delta, \tH_\delta \big)= \big(G* \oE_\delta, G*\oH_\delta \big)  \mbox{ in } B_{r_3}. 
\end{equation}
Since $G\circ F (x) = \big(r_3^2/ r_2^2 \big) x$ and  $G_*F_*  = (G \circ F)_*$,  it follows from \eqref{def-T*} and \eqref{choice-Br1} that 
\begin{equation}\label{important-identity}
(G_*F_*\eps_\delta, G_*F_*\mu_\delta) =  (G_*F_*\eps_O, G_*F_*\mu_O) = (I, I) \mbox{ in } B_{r_3}. 
\end{equation}
Set 
\begin{equation}\label{def-Data}
Data(j, \delta): =  \left(\frac{1}{\delta} \| (E_\delta, H_\delta)\|_{L^2(B_{R_0} \setminus B_{r_3})} \| j\|_{L^2} + \| j\|_{L^2}^2 \right)^{1/2}. 
\end{equation}
Applying Lemma~\ref{lem-stability} to $D = B_{r_2} \setminus B_{r_1}$, we have 
\begin{equation}\label{est-F}
\| (E_\delta, H_\delta)\|_{[L^2(B_{R_0})]^2}^2 \le C Data(j, \delta)^2.  
\end{equation}
Here and in what follows in the proof, $C$ denotes a positive constant independent of $\delta$ and  $j$ and we assume that $\ell > 10$. 

\medskip
The proof now 
is divided into two steps. 
\begin{itemize}
\item Step 1:  We prove that if $\ell$ is large enough then 
\begin{equation}\label{step1}
\| (E^{(1)}_\delta - E_\delta) \times \nu, (H^{(1)}_\delta - H_\delta) \times \nu\|_{H^{-1/2}(\dive_\Gamma, \partial B_{2 r_2})}  \le C \delta^{\gamma + 1/2} Data(j, \delta). 
\end{equation}

\item Step 2: Define 
\begin{equation*}
(\cE_\delta, \cH_\delta)  = \left\{ \begin{array}{cl} (E_\delta, H_\delta) & \mbox{ in } \mR^3 \setminus B_{r_3}, \\[6pt]
(E_\delta, H_\delta) - \big(\oE_\delta - \tE_\delta,  \oH_\delta - \tH_\delta \big)  & \mbox{ in } B_{r_3} \setminus B_{2 r_2}, \\[6pt]
\big(\tE_\delta, \tH_\delta \big) & \mbox{ in } B_{2 r_2}.  
\end{array}\right.
\end{equation*}
We prove that if \eqref{step1} holds  then 
\begin{equation}
\| (\cE_\delta, \cH_\delta) -  (E, H) \|_{L^2(B_{R} \setminus B_{r_3})} \le C \delta^{\gamma} \| j \|_{L^2}. 
\end{equation}
\end{itemize}

It is clear that the conclusion follows after Step 2. 

\medskip
\noindent \underline{Step 1:} Using the fact that  
$$
i k F^{-1}_* \tmu_O + i \delta I = i k \big(F^{-1}_*\tmu_O + (\delta/k) F^{-1}_* F_* I   \big)  \mbox{ in } B_{r_2} \setminus B_{r_1}
$$
and
$$
- i k F^{-1}_* \teps_O + i \delta I = -  i k \big(F^{-1}_*\teps_O - (\delta/k) F^{-1}_* F_* I   \big)  \mbox{ in } B_{r_2} \setminus B_{r_1}, 
$$  
and applying Lemma~\ref{lem-CV}, we have 
\begin{equation}\label{sys-EH1}
\left\{\begin{array}{cl}
\nabla \times \oE_\delta = i k \tmu_O \oH_\delta + i \delta F_* I  \oH_\delta  &   \mbox{ in } B_{r_3} \setminus B_{r_2} , \\[6pt]
\nabla \times \oH_\delta = -i k   \teps_O \oE_\delta +  i \delta F_* I \oE_\delta & \mbox{ in } B_{r_3} \setminus B_{r_2},  
\end{array} \right.  \quad 
\end{equation}
and 
\begin{equation}\label{transmission-r2}
\big(\oE_\delta \times \nu, \oH_\delta \times \nu \big)= \big(E_\delta \times \nu, H_\delta \times \nu \big) \Big|_{\rm ext} \mbox{ on } \partial B_{r_2}. 
\end{equation}
In \eqref{transmission-r2},  we use the fact that  $F(x) = x$ on $\partial B_{r_2}$. Set  
\begin{equation}\label{def-eps-mu}
(\eps, \mu) = \left\{\begin{array}{cl} (\teps_O, \tmu_O) & \mbox{ in } B_{r_3} \setminus B_{r_2}, \\[6pt]
(I, I) & \mbox{ otherwise.}
\end{array} \right.
\end{equation}
Let $\big(\obE_\delta, \obH_\delta \big) \in [H_{\loc}(\mR^3)]^2$  be the unique outgoing solution to
\begin{equation}\label{def-EH1}
\left\{\begin{array}{cl}
\nabla \times \obE_\delta = i k \mu \obH_\delta + i \delta  \mathds{1}_{B_{r_3} \setminus B_{r_2}} F_*I \oH_\delta& \mbox{ in } \mR^3, \\[6pt]
\nabla \times  \obH_\delta = - i k \eps \obE_\delta + i  \delta  \mathds{1}_{B_{r_3} \setminus B_{r_2}} F_*I \oE_\delta & \mbox{ in } \mR^3.
\end{array}\right.
\end{equation}
Here $\mathds{1}_{D}$ denotes the characteristic function of a subset $D$ of $\mR^3$. Applying  Lemma~\ref{lem-stability} to $(\obE_\delta, \obH_\delta)$ with $f = i \delta  \mathds{1}_{B_{r_3} \setminus B_{r_2}} F_*I \oH_\delta$, $g = i  \delta  \mathds{1}_{B_{r_3} \setminus B_{r_2}} F_*I \oE_\delta $, and $h_1 = h_2 = 0$, and using \eqref{est-F},   we obtain 
\begin{equation}\label{part1-1}
\| \big( \obE_\delta, \obH_\delta \big) \|_{ H(\curl, B_{r_3} \setminus B_{r_2})} \le C \delta Data(j, \delta). 
\end{equation}
Set 
\begin{equation*}
(\ocE_\delta, \ocH_\delta) = 
\left\{\begin{array}{cl}
\Big(\oE_{\delta} - E_\delta - \obE_\delta, \oE_{\delta} - H_\delta - \obH_\delta \Big) & \mbox{ in } B_{r_3}  \setminus B_{r_2}, \\[6pt]
\Big(-\obE_{\delta}, - \obH_\delta \Big)  &  \mbox{ in } B_{r_2}. 
\end{array}\right.
\end{equation*}
It follows from \eqref{sys-EH1} and \eqref{def-EH1} that 
\begin{equation}\label{ocE}
\left\{\begin{array}{cl}
\nabla \times \ocE_\delta = i k \mu \ocH_\delta & \mbox{ in } B_{r_3} \\[6pt]
\nabla \times \ocH_\delta = - i k\eps \ocE_\delta & \mbox{ in } B_{r_3}. 
\end{array}\right.
\end{equation}
Applying Lemma~\ref{lem-elliptic}, we have, for $1 \le  a \le 3$,  
\begin{equation*}
\dive \big(\eps \nabla \ocE_{\delta, a} \big) = -  \dive \big(\partial_a \eps  \ocE_{\delta} + i k \eps \epsilon^a \mu  \ocH_{\delta} \big) \quad \mbox{ in } B_{r_3}
\end{equation*}
and 
\begin{equation*}
 \dive \big(\mu \nabla \ocH_{\delta, a} \big)   = -  \dive \big(\partial_a \mu  \ocH_{\delta} - i k \mu \epsilon^a \eps_{\delta} \ocE_\delta \big) \quad \mbox{ in } B_{r_3}, 
\end{equation*}
where $\epsilon^a_{bc}$ ($1 \le a, \, b, \, c \le 3$) denote the usual Levi Civita permutation, see \eqref{Levi-Civita}. 
Let $q$ be the constant in Lemma~\ref{lem-3S} with $m = 3$, $n=6$,
$M^1 = M^2 = M^3 = \eps$, and $M^4 = M^5 = M^6 = \mu$. Define, for $0 < r \le r_3$, 
$$
\| \ocE \|_{\bH(\pB_r)} = \| \ocE \|_{H^{1/2}(\pB_r)} + \| \eps \nabla \ocE \cdot e_r \|_{H^{-1/2}(\pB_r)}, 
$$
and
$$
\| \ocH \|_{\bH(\pB_r)} = \| \ocH \|_{H^{1/2}(\pB_r)} + \| \mu \nabla \ocH \cdot e_r \|_{H^{-1/2}(\pB_r)}. 
$$
By Lemma~\ref{lem-3S}, there exists some positive constant $C$ independent of $\delta$ such that 
\begin{equation}\label{part1-2-*}
\| \big(\ocE_\delta, \ocH_\delta \big) \|_{\bH(\partial B_{2r_2}) }  \le C \| \big(\ocE_\delta, \ocH_\delta \big)\|_{\bH(\partial B_{r_2/2}) }^\alpha \| \big(\ocE_\delta, \ocH_\delta \big)\|_{ \bH(\partial B_{4 r_2}) }^{1 - \alpha},  
\end{equation}
with 
\begin{equation}\label{def-alpha-part1}
\alpha = \frac{(2r_2)^{-q} - (4 r_2)^{-q}}{(r_2/2)^{-q} - (4 r_2)^{-q}} =  \frac{2^{-q} - 4 ^{-q}}{2^{q} - 4^{-q}}. 
\end{equation}
Since $\eps = \mu = I$ in $B_{r_2} \cup ( B_{r_3} \setminus B_{r_3/4})$ (recall that $\ell > 10$), it follows from  \eqref{ocE} that   
\begin{equation}\label{Laplace}
\Delta \ocE_\delta + k^2 \ocE_\delta = \Delta \ocH_\delta + k^2 \ocH_\delta  = 0  \mbox{ in } B_{r_2} \cup ( B_{r_3} \setminus B_{2 r_2}). 
\end{equation}
Applying Lemma~\ref{lem-3S1}, we have 
\begin{equation}\label{part1-2-1}
\| (\ocE, \ocH) \|_{\bH(\partial B_{4r_2})} \le C \| (\ocE, \ocH) \|_{\bH(\partial B_{2r_2})}^{\beta} \| (\ocE, \ocH) \|_{\bH(\partial B_{r_3/2})}^{1 - \beta}, 
\end{equation}
where 
\begin{equation}\label{def-beta}
\beta = \ln \Big( \frac{r_3}{4 r_2} \Big) \Big/ \ln \Big( \frac{r_3}{2 r_2} \Big). 
\end{equation}
Combining \eqref{part1-2-*} and \eqref{part1-2-1} yields 
\begin{equation}\label{part1-2}
\| \big(\ocE_\delta, \ocH_\delta \big) \|_{\bH(\partial B_{2r_2}) }  \le C \| \big(\ocE_\delta, \ocH_\delta \big)\|_{\bH(\partial B_{r_2/2}) }^\rho \| \big(\ocE_\delta, \ocH_\delta \big)\|_{ \bH(\partial B_{r_3/2}) }^{1 - \rho},   
\end{equation}
with 
\begin{equation}\label{def-rho}
\rho = \frac{\alpha}{1 - (1 - \alpha) \beta }. 
\end{equation}

On the other hand, from \eqref{Laplace}, we derive that 
\begin{equation}\label{part1-L2}
\| (\ocE, \ocH) \|_{\bH(\pB_{r_2/2})} \le C \| (\ocE, \ocH) \|_{L^2(B_{r_2})}  \mbox{  and  }  \| (\ocE, \ocH) \|_{\bH(\pB_{r_3/2})} \le C \| (\ocE, \ocH) \|_{L^2(B_{r_3} \setminus B_{r_3/4})}. 
\end{equation}
From   \eqref{est-F}, \eqref{part1-1}, \eqref{part1-2}, and \eqref{part1-L2}, we obtain 
\begin{equation}\label{part1-3}
\| \big(\ocE_\delta, \ocH_\delta \big)\|_{\bH(\partial B_{2r_2})} \le  C \delta^{\rho}  Data(j, \delta). 
\end{equation}
By taking $l$ large enough, we derive from  \eqref{def-alpha-part1}, \eqref{def-beta}, and \eqref{def-rho} that $\rho > 1/2 + \gamma$ if $r_3 > l r_2$. The conclusion of Step 1 follows.

\medskip 

\noindent \underline{Step 2:} We have, since $G(x) = x$ on $\partial B_{r_3}$,  
\begin{equation*}
[\cE_\delta \times \nu]  = (\oE_\delta - \tE_\delta) \times \nu  = 0  \mbox{ on } \partial B_{r_3}, 
\end{equation*}
and 
\begin{equation*}
[\cH_\delta \times \nu]  = (\oH_\delta - \tH_\delta) \times \nu  = 0  \mbox{ on } \partial B_{r_3}. 
\end{equation*}
Applying Lemma~\ref{lem-CV},  we obtain 
\begin{equation}\label{sys-cEcH}
\left\{\begin{array}{cl}
\nabla \times  \cE_\delta = -  i k  \cH_\delta & \mbox{ in } \mR^3 \setminus \partial B_{2 r_2}, \\[6pt]
\nabla \times  \cH_\delta = i k \cE_\delta + j & \mbox{ in } \mR^3 \setminus \partial B_{2 r_2}.  
\end{array}\right. 
\end{equation}
Applying Lemma~\ref{lem-stability},  we derive  from \eqref{step1} that 
\begin{equation}\label{p1}
\| (\cE_\delta, \cH_\delta)\|_{H(\curl, B_{R_0} \setminus \partial B_{r_3})} \le C \Big(\| j\|_{L^2} + \delta^{\gamma + 1/2}  Data(j, \delta) \Big).
\end{equation}
Since $\gamma > 0$ and $(\cE_\delta, \cH_\delta) = (E_\delta, H_\delta)$ in $\mR^3 \setminus B_{r_3}$, it follows from \eqref{def-Data} and \eqref{p1} that 
\begin{equation*}
\| (\cE_\delta, \cH_\delta) \|_{H(\curl, B_{R_0} \setminus B_{r_3})]^2} \le C \| j\|_{L^2}. 
\end{equation*}
We derive from \eqref{def-Data} that 
\begin{equation}
Data(j, \delta) \le C \delta^{-1/2} \| j\|_{L^2};  
\end{equation}
which in turn implies, by \eqref{step1}, 
\begin{equation}\label{est-L2-part2}
\| (E^{(1)}_\delta - E_\delta) \times \nu, (H^{(1)}_\delta - H_\delta) \times \nu\|_{\big[H^{-1/2}(\dive_\Gamma, \partial B_{2 r_2}) \big]^2 }  \le C \delta^{\gamma} \| j\|_{L^2}. 
\end{equation}
It is clear from \eqref{eq-EH} and \eqref{sys-cEcH} that  
\begin{equation*}
\left\{\begin{array}{cl}
\nabla \times  (\cE_\delta  - E) = -  i k  (\cH_\delta - H) & \mbox{ in } \mR^3 \setminus \partial B_{2 r_2}, \\[6pt]
\nabla \times  (\cH_\delta - H) = i k (\cE_\delta - E) & \mbox{ in } \mR^3 \setminus \partial B_{2 r_2}.  
\end{array}\right. 
\end{equation*}
Using \eqref{est-L2-part2} and  applying Lemma~\ref{lem-stability} for the system $(\cE_\delta - E, \cH_\delta - H)$, one obtains the conclusion of Step 2. 

\medskip 
The proof is complete. \proofend

\begin{remark}  \fontfamily{m} \selectfont  The definition of $(\cE_\delta, \cH_\delta)$ is one of the key points of the proof. The idea is to remove from $(E_\delta, H_\delta)$ the term $(\oE_\delta - \tE_\delta, \oH_\delta - \tH_\delta)$ in $B_{r_3} \setminus B_{2r_2}$; which is singular in general.  This is the spirit of the removing of localized singularity technique introduced in \cite{Ng-Negative-Cloaking, Ng-Superlensing}.
\end{remark}

\section{Further discussion} \label{sect-further-discussion}

 In Theorem~\ref{thm1}, the requirement that $F$ is the Kelvin transform with respect to $\partial B_{r_2}$ can be relaxed. In fact, as seen in the proof of Theorem~\ref{thm1}, one can replace the Kelvin transform by any transformation $F: B_{r_2} \setminus \bar B_{r_1} \to B_{r_3} \setminus \bar B_{r_2}$ such that $i)$ $F(x)  = x$ on $\partial B_{r_2}$; $ii)$ There exists a diffeomorphism extension of $F$, which is still denoted by  $F$, from $B_{r_2} \setminus \{0\}$  onto $\mR^3 \setminus \bar B_{r_2}$;  $iii)$ There exists a diffeomorphism $G: \mR^3 \setminus \bar B_{r_3} \to B_{r_3} \setminus \{0 \}$  such that  $G \in C^1(\mR^3 \setminus B_{r_3})$,  $G(x) = x \mbox{ on } \partial B_{r_3}$, and $G \circ F : B_{r_1}  \to B_{r_3} \mbox{ is a diffeomorphism if one sets } G\circ F(0) = 0$. In this context, the first layer in $B_{r_2} \setminus B_{r_1}$ is also given by \eqref{first-layer} and the second layer in $B_{r_1}$ is changed correspondingly  by
\begin{equation}\label{second-layer}
\big( F^{-1}_*G^{-1}_*I, F^{-1}_*G^{-1}_*I\big). 
\end{equation}
Set 
\begin{equation}\label{def-eps-mu-1}
(\eps_\delta, \mu_\delta) = \left\{ \begin{array}{cl} 
\big( \teps_O , \tmu_O \big)& \mbox{ in } B_{r_3} \setminus B_{r_2}, \\[6pt]
\big( F^{-1}_* \teps_O +  i  \delta I, F^{-1}_* \tmu_O  + i\delta I  \big) & \mbox{ in } B_{r_2} \setminus B_{r_1}, \\[6pt]
\big( F^{-1}_*G^{-1}_*I, F^{-1}_*G^{-1}_*I\big) & \mbox{ in } B_{r_1},\\[6pt]
 \big(I,  I \big) & \mbox{ in } \mR^3 \setminus B_{r_3}.  
\end{array} \right. 
\end{equation} 

We have 
\begin{proposition}\label{pro1} Let $R_0>0$,   $j \in \big[L^{2}(\mR^3) \big]^3$ with $\supp j \subset \subset B_{R_0} \setminus  B_{r_{3}}$ and let  $(E_\delta, H_\delta) \in  \big[H_{\loc}(\curl, \mR^3) \big]^2$ be the unique  outgoing solution to \eqref{eq-EHdelta} where $(\eps_\delta, \mu_\delta)$ is given by \eqref{def-eps-mu-1} and let $(E, H) \in \big[H_{\loc}(\curl, \mR^3) \big]^2$ be the unique  outgoing solution to \eqref{eq-EH}. Given $0 < \gamma < 1/2$, there exists  a positive constant  $\ell = \ell(\gamma)> 0$, depending {\bf only} on  the elliptic  constant of $\teps_O$ and $\tmu_O$ in $B_{2 r_2} \setminus B_{r_2}$  and  $\| (\teps_O, \tmu_O) \|_{W^{2, \infty}(B_{4 r_2})}$ such that if $r_3 > \ell r_2$ then   
\begin{equation}\label{key-point}
\| (E_{\delta}, H_\delta) - (E, H) \|_{H(\curl, B_R \setminus B_{r_3})} \le C_R \delta^\gamma \| j\|_{L^2}, 
\end{equation}
for some positive constant $C_R$ independent of $j$ and $\delta$. 
\end{proposition}

The constants $\gamma$ and $\ell(\gamma)$ in Proposition~\ref{pro1} can be chosen as the ones in Theorem~\ref{thm1}. The proof of Proposition~\ref{pro1} follows the same line as the one of Theorem~\ref{thm1} and is omitted.

\providecommand{\bysame}{\leavevmode\hbox to3em{\hrulefill}\thinspace}
\providecommand{\MR}{\relax\ifhmode\unskip\space\fi MR }
\providecommand{\MRhref}[2]{%
  \href{http://www.ams.org/mathscinet-getitem?mr=#1}{#2}
}
\providecommand{\href}[2]{#2}

%\bibliographystyle{amsplain}
%\bibliography{/Dropbox/Bib/bib1}
%\bibliography{/Users/hoaiminhnguyen/Dropbox/Bib/bib1}
%\bibliography{/Users/hnguyen/Dropbox/Bib/bib1}
%\bibliography{mybib}

\end{document}